\documentclass[12pt, a4paper, twoside]{amsart}
\usepackage{amsmath, amssymb, amsthm, amsfonts}
\usepackage[a4paper, left=2.7cm, right=3cm, top=3cm, bottom=3cm]{geometry}

\newcommand{\R}{\mathbb{R}} 
\newcommand{\C}{\mathbb{C}}
\newcommand{\Z}{\mathbb{Z}}
\newcommand{\define}{\mathrel{\mathop:}=}

\newcommand{\App}{\mathcal{A}}
\newcommand{\sW}{W} 
\newcommand{\aW}{{W^a}} 
\newcommand{\Cf}{\mathcal{{C}}_{f}} 
\newcommand{\Cfm}{\mathcal{{C}}_{f}^{\mathrm{op}}} 
\newcommand{\fa}{{\bf{c}_f}} 
\newcommand{\Cp}{\mathcal{{C}}_{p}} 
\newcommand{\Cpm}{\mathcal{{C}}_{p}^{\mathrm{op}}} 
\newcommand{\Hd}{\overline{H}}
\newcommand{\Ghat}{\hat{\mathcal{G}}} 
\newcommand{\G}{\mathcal{G}} 
\newcommand{\dconv}{\mathrm{{conv}}} 
\newcommand{\wgt}{\mathrm{{wgt}}} 
\newcommand{\fw}{{\omega}} 
\newcommand{\type}{\mathrm{type}} 

\newtheorem{thm}{Theorem}[section]
\newtheorem*{thm*}{Theorem}
\newtheorem{prop}[thm]{Proposition}
\newtheorem*{prop*}{Proposition}
\newtheorem{lemma}[thm]{Lemma}
\newtheorem*{lemma*}{Lemma}
\newtheorem{ex}[thm]{Example}

\newtheorem{cor*}[thm]{Corollary}
\newtheorem{definition}[thm]{Definition}
\newtheorem*{definition*}{Definition}
\newtheorem{rem}[thm]{Remark}
\newenvironment{prf}{\begin{proof}[Proof]}{\end{proof}}

\newenvironment{proofMainThm}{\begin{proof}[Proof of \ref{mainthm}]}{\end{proof}}
\newenvironment{proofPropA}{\begin{proof}[Proof of \ref{propA}]}{\end{proof}}
\newenvironment{proofPropB}{\begin{proof}[Proof of \ref{propB}]}{\end{proof}}

\newtheorem*{propA}{Proposition \ref{propA}}
\newtheorem*{propB}{Proposition \ref{propB}}
\newtheorem*{mainthm}{Theorem \ref{mainthm}}

\begin{document}

\title{Kostant convexity for affine buildings}
\author{Petra Hitzelberger\\ February 2008 }

\begin{abstract}
We prove an analogue of Kostants convexity theorem for thick affine buildings and give an application for groups with affine $BN$-pair. Recall that there are two natural retractions of the affine building onto a fixed apartment $A$: The retraction $r$ centered at an alcove in $A$ and the retraction $\rho$ centered at a chamber in the spherical building at infinity. We prove that for each special vertex $x\in A$ the set $\rho(r^{-1}(W.x))$ is a certain convex hull of $W.x$. The proof can be reduced to a statement about Coxeter complexes and heavily relies on a character formula for highest weight representations of algebraic groups.
\end{abstract}

\maketitle

\section{Introduction}

Kostant \cite{Kostant} proved a convexity theorem for symmetric spaces, generalizing a well known theorem of Schur \cite{Schur}. Let $G\!/\!K$ be a symmetric space and $\mathfrak{g}=\mathfrak{k}\oplus\mathfrak{p}$ the Cartan decomposition of the associated Lie algebra. Let $T$ be a maximal flat of $G\!/\!K$ and $\sW=N_K(T)\!/\!{Z_K(T)}$ the spherical Weyl group of $G\!/\!K$ which acts on $T$. Let $\pi$ denote the Iwasawa projection onto $T$, that is $\pi: ank \mapsto a$. Kostant's result is
\begin{thm*}[\cite{Kostant}]
Let $K.x$ be the $K$ orbit of $x\in G\!/\!K$. Then $$\pi(K.x) = \dconv(W.x).$$
\end{thm*}

Affine buildings are in many ways similar to symmetric spaces. Often results for symmetric spaces have an analogue in the theory of affine buildings. Some results can be stated simultanously. Therefore it is a natural question whether a convexity result in the spirit of Kostant's can be found for arbitrary affine buildings $X$. It turns out that an analog statement is true, assuming thickness of $X$. The flat $T$ is replaced by an apartment $A$ (together with a fixed origin $0$), which is a flat subspace of the affine building $X$, equipped with a natural action of the spherical Weyl group $\sW$. The role of the Iwasawa projection is played by a retraction of $X$ onto the apartment $A$. If the automorphism group $\mathrm{Aut}(X)$ is big enough the set $K.x$ is replaced by the orbit of $x$ under the stabilizer of $0$ in $\mathrm{Aut}(X)$.

\subsection*{The main result} We fix an apartment $A$ in $X$. 
Let $S$ in $A$ be a Weyl chamber, that is a translate of the metric closure of a connected component of $A\setminus \{H_{i,0}, i\in I\}$ under the affine Weyl group $\aW$. The parallel class $\partial S$ of the Weyl chamber $S$ is a chamber in the spherical building $\Delta=\partial X$ at infinity. Let $\Cf$ denote the fundamental Weyl chamber. Translates of a fundamental domain of the action of $\aW$ are alcoves in $A$ and we will write $\fa$ for the fundamental one.

There are the following two retractions onto $A$: The canonical retraction $r$ centered at the fundamental alcove $\fa\in A$ and the retraction $\rho$ centered at $\partial\Cfm \in \partial A$. See page~\pageref{retractions} for definitions.

Denote by $\fw_i$ the fundamental co-weights. For all $i\in I$ and $k\in \Z$ the set $\Hd_{i,k}\define\{ x\in A\vert \langle x^\vee, \fw_i\rangle=k \}$ is called a dual hyperplane. Usually $W$-convex sets in a building are defined as intersections of finitely many half apartments. Analogously we define a convex set\footnote{One could also use the somewhat redundant name \emph{dual convexity}.} to be an intersection of finitely many dual half apartments. The convex hull of a set $Y$, denoted by $\dconv(Y)$, is the smallest convex set containing~$Y$. Our main result reads as follows

\begin{mainthm}
Let $X$ be a thick affine building and $A$ an apartment in $X$. Let $r$ and $\rho$ be defined as above. For all special vertices $x$ in $A$, 
$$\rho(r^{-1}(\sW .x))\;=\dconv(W.x)\cap (x+Q),$$
where $Q=Q(R^\vee)$ the fundamental co-weight lattice of the underlying root system $R$.
\end{mainthm}

By extending the maps $\rho$ and $r$ to galleries in $X$ starting at $0$, the proof of the theorem can be reduced to the following two facts:
\begin{propA}
All vertices in $\dconv(W.x)$ of the same type as $x$ are endpoints of positively folded galleries of fixed type $t$, where $t$ is the type of a minimal gallery connecting the origin $0$ with $x$. All such endpoints are contained in this set.
\end{propA}

\begin{propB}
Every positively folded gallery in the proposition above has a pre-image under $\rho$ which is a minimal gallery starting at $0$. 
\end{propB}

In \cite{GaussentLittelmann} Gaussent and Littelmann study irreducible highest weight representations of semisimple complex or compact real Lie groups. They establish a character formula in terms of positively folded galleries. Proposition~\ref{propA} is an easy consequence of their character formula. Note that the first proposition is just a statement about affine Coxeter complexes. One is therefore completely independent of the ambient building. Our main result holds for \emph{all} thick affine buildings, which need not be associated to an algebraic group, even though the proof relies on deep results about representations of algebraic groups.

\subsection*{An application}
Let $X, A$,and  $\Cfm$ be as above. Assume $G$ is a group acting on $X$ by automorphisms. Assume further that $G$ has a split spherical and an affine $BN$-pair. Hence the spherical $B$ splits as $B=UT$. Let $K$ be the stabilizer of $0$ in $G$. Then $G$ has an Iwasawa decomposition
$ G=UTK $, where $U$ is the unipotent radical and $T$ is the group of translations in $A$.  Therefore we can identify special vertices in $X$ of the same type as $0$ with cosets $utK$ in $G$, where $u\in U$ and $t\in T$. The origin $0$ corresponds to $K$.  It is easy to verify that $\rho$ is the map $utK\mapsto tK$. The set $r^{-1}(W.x)$ is the same as the $K$-orbit of $x$, for all $x\in A$. Note that the vertices in $A$ having the same type as $0$ are of the form $tK$ for some $t\in T$. The following corollary is just a reformulation of our main result. 

\begin{thm*}
Given $x=tK\in X$ the following is true: 
$$\rho(KtK)=\dconv(W.tK).$$
Or equivalently, since $\rho^{-1}(t'K)=Ut'K$,
$$
\emptyset \neq Ut'K \cap KtK \Leftrightarrow t'K \in \dconv(W. tK).
$$
\end{thm*}

The proof of "$\Rightarrow$" is well known and can be found in \cite{BruhatTits}. For partial results on "$\Leftarrow$" and related questions see for example \cite{GaussentLittelmann, MirkovicVilonen}. Denote by $M$ the matrix describing the Satake isomorphism in terms of natural bases of source and target. Sometimes the implication from the right to the left is stated in terms of positivity of the coefficients of $M$. Rapoport proved such a positivity result in \cite{Rapoport}, which was generalized by Schwer in \cite{Schwer}. Compare also Silberger, \cite{Silberger}, where a proof of the equivalence $\emptyset\neq Ut'K \cap KtK \Leftrightarrow t'K \in \dconv(tK)$ for algebraic groups defined over $p$-adic fields was announced.

\subsection*{The paper is organized as follows}
We set notations and state the needed definition and results in Section 2. Here we also define positively folded galleries, convexity and state the character formula mentioned above. Section 3 is devoted to the proof of our main result. In the last Section the application is discussed.

\thanks{The author would like to thank L. Kramer and Christoph Schwer for many helpful diskussions and remarks.

This work is part of the author's doctoral thesis, written under the supervision of L. Kramer and R. Weiss at the University of M\"unster, Germany. During the preparation of this paper the author was supported by the \emph{Studienstiftung des deutschen Volkes} and the \emph{DFG}.

\section{Preliminaries, definitions and notation}

\subsection*{Alcoves and Weyl chambers}

For the following definitions see for example \cite{Bourbaki4-6, Brown}. Let $V$ be a vector space over $\R$ and denote by $V^*$ its dual. We denote by $\langle\alpha, \beta^\vee\rangle$ the evaluation of $\beta\in V^*$ on $\alpha\in V$. Recall the following 
\begin{definition}
A subset $R$ of $V$ is a \emph{root system} in $V$ if the following holds
\begin{itemize}
\item $R$ is finite, $0\notin R$ and $R$ generates $V$
\item for all $\alpha\in R$ exists $\alpha^\vee\in V^*$ such that $\langle \alpha, \alpha^\vee\rangle = 2$ 
\item the reflection $r_{\alpha, \alpha^\vee}: x\mapsto x-\langle x, \alpha^\vee\rangle \alpha$ leaves $R$ invariant
\item for all $\alpha\in R$ the set $\alpha^\vee(R)$ is contained in $\Z$. 
\end{itemize}
\end{definition}
Note that $R^\vee=\{\alpha^\vee : \alpha\in R\}$ is a root system in $V^*$. We refer to $R^\vee$ as the \emph{inverse root system}. If $B$ is a basis of $R$ then $B^\vee=\{\alpha^\vee\}_{\alpha\in B}$ is a basis of $R^\vee$. The elements of the dual basis $\{\fw_\alpha\}_{\alpha\in B}$ of $B^\vee$ in $V$ are called \emph{fundamental co-weights}.

Let $X$ be an (irreducible, simplicial) affine building with apartment system $\App$. Let $R$ denote the associated irreducible root system contained in $V$. Let $B=\{\alpha_i, i\in I\}$ be a root basis of $R$. Further let $\aW$ denote the \emph{affine} and $\sW$ the \emph{spherical Weyl group}.

Let $E$ denote the affine space underlying the vector space $V$ on which $R$ is defined. The \emph{hyperplane} associated to $\alpha_i\in B$ is $H_{i,0}=\{x\in E : \langle x, \alpha^\vee_i \rangle =0\}$. By \emph{Weyl chambers} we mean translates of the closures of connected components of $E\!\setminus\! \{H_{i,0}, i\in I\}$ under the affine Weyl group $\aW$. 

A hyperplane is an image of some $H_{i,0}$ under $\aW$. For each hyperplane $H$ there exist $\alpha\in R$ and $k\in\Z$ such that $H=H_{\alpha,k}\define\{x\in E \;\vert\; \langle x, \alpha^\vee\rangle = k \}$. We will refer to closures of connected components of $E\setminus\{H_{\alpha, k}, \alpha\in R, k\in\Z\}$ as \emph{alcoves}. Denote by $\alpha_0$ the highest root in $R$. Then 
$$
\fa=\{x\in E \;\vert\; \langle x,\alpha_0^\vee\rangle \leq 1 \text{ and } \langle x, \alpha_i^\vee\rangle \geq 0 \text{ for all } i\in I \}
$$
is a \emph{fundamental alcove}.
The parallel class $\partial S$ of the Weyl chamber $S$ is an alcove in the \emph{spherical building $\Delta=\partial X$ at infinity}. 
A \emph{half appartment} $H_{\alpha,k}^{\pm}$ is one of the two half spaces $\{x\in E : \langle x, \alpha^\vee\rangle \geq k\}$ and $\{x\in E : \langle x,\alpha^\vee\rangle \leq k\}$.
Let $S$ be a Weyl chamber. The intersection of all faces of $S$ is the \emph{vertex of $S$}. The set of Weyl chambers in $E$ sharing the same vertex is in one-to-one correspondence with the chambers of the spherical Coxeter complex associated to $\sW$.

Denote by $\Cf$ the \emph{fundamental Weyl chamber}. By definition
\begin{align*}
\Cf 	&= \{x\in E\;\vert\; \langle x,\alpha_i^\vee\rangle \geq 0 \;\text{ for all } i\in I\}\\
	&= \{x\in E\;\vert\; x=\sum_{i\in I} k_i\fw_i \;\text{ with } k_i\geq 0 \;\text{ for all } i\in I\}.
\end{align*}

\subsection*{Cones and dual cones}
The \emph{positive cone $\Cp$}, which is the dual cone of $\Cf$, will play an important role in the definition of the convexity we use in the formulation of our main result.

The \emph{positive cone $\Cp$} can be described as follows:
\begin{align*}
\Cp 	&= \{x\in E\;\vert\; \langle \fw_i, x^\vee\rangle \geq 0 \;\text{ for all } i\in I\}\\
	&= \{x\in E\;\vert\; x=\sum_{i\in I} k_i\alpha_i \;\text{ with } k_i\geq 0 \;\text{ for all } i\in I\}.
\end{align*}


\subsection*{Retractions}\label{retractions}

Fix an apartment $A$ in $X$ and a special vertex $0$ in $A$.
There are two natural retractions from $X$ onto $A$, which can be characterized as follows: 
\begin{definition}
Let $c$ be an alcove in $A$. The \emph{canonical retraction onto $A$ centered at $c$}, denoted by $r_{A,c}:X\rightarrow A$, is characterized by the following two properties:
\begin{enumerate}
\item The alcove $c$ is fixed pointwise by $r_{A,c}$ and each apartment containing $c$ is mapped isomorphically onto $A$.
\item The map $r_{A,c}$ is distance diminishing on $X$.
\end{enumerate}
\end{definition}

\begin{definition}
Let $S$ be a Weyl chamber in $A$. The \emph{retraction onto $A$ centered at $\partial S$}, which will be denoted by $\rho_{A,S}:X\rightarrow A$, is the unique map such that the following conditions are satisfied:
\begin{enumerate}
\item For each apartment $B$ containing a sub-Weyl chamber of $S$ the restriction of $\rho_{A,S}$ to $B$ is an isomorphism onto $A$ fixing $A\cap B$ pointwise.
\item The map $\rho_{A,S}$ diminishes distances in $X$.
\end{enumerate}
\end{definition}

One has the following connection between these two retractions:

\begin{thm*}[\cite{Brown}, p.171]
Let $\rho_{A,S}$ be the retraction onto $A$ centered at $\partial S \in A$. Given an alcove $d\in X$ there exists an alcove $c\in A$ such that $$\rho_{A, S}(d)=r_{A,c}(d).$$
\end{thm*}

We identify $\sW$ with $\mathrm{Stab}_\aW(0)$. The Weyl chamber opposite $\Cf$ in $A$ is denoted by $\Cfm$. For simplicity of notation, we write $r$ instead of $r_{A,\fa}$ and let $\rho$ stand for $\rho_{A,\Cfm}$. 

Choosing apartments $A_1, A_2$ such that $\partial\Cfm\subset\partial A_i$  for $i=1$ and $2$, the corresponding retractions $\rho_{A_1, \Cfm}$ and $\rho_{A_2, \Cfm}$ are compatible in the sense of the following lemma.
\begin{lemma}\label{compatibleRetractions}
Given apartments $A_i, i=1,\ldots, n$, such that $\partial\Cfm\subset \partial A_i$ for all $i$. Denote by $\rho_i$ the retraction $\rho_{A_i, \Cfm}$. Then 
$$
\rho_1\circ\rho_2\circ\ldots\circ\rho_n = \rho_1 .
$$
\end{lemma}
\begin{proof}
For all $i$ the retraction $\rho_i$ maps apartments containing representatives of $\partial\Cfm$ isomorphically onto $A_i$. As a set $X$ equals the union of all such apartments. Every such apartment is isomorphically mapped onto the apartment $A_1$ by $\widetilde{\rho}\define \rho_1\circ\rho_2\circ\ldots\circ\rho_n$. A second property of the map $\widetilde{\rho}$ is that it is distance diminishing, since each $\rho_i$ is. Therefore $\widetilde{\rho}:X\mapsto A_1$ is the unique map defined by these two properties, and therefore equals $\rho_1$.
\end{proof}

For further details see for example \cite{Brown}.


\subsection*{Positively folded galleries}

In order to have more control over stammering galleries we introduce the slightly more general notion of combinatorial galleries. One can think of a combinatorial gallery as a gallery in the usual sense together with a fixed type.

\begin{definition}\label{galleryDef}
A \emph{combinatorial gallery} $\gamma$ in $X$ is a sequence of faces and alcoves
$$\gamma = (c_0', c_0, c_1', c_1, c_2', \ldots ,c_n', c_n, c_{n+1}' )$$ 
such that 
\begin{itemize}
\item the \emph{source} $c_0'$ and the \emph{target} (or \emph{endpoint}) $c_{n+1}'$ are vertices in $X$,
\item $c_i$ is an alcove, for all $i=0\ldots n$, and
\item $c_i'$ is a codimension one face of $c_{i-1}$ and $c_i$ for all $i=1\ldots n-1$.
\end{itemize}
The \emph{type of a combinatorial gallery} $\gamma$ is the list $t=(t_0,t_1,\ldots, t_n)$ of types $t_i$ of the faces $c_i'$.
We say $\gamma$ is \emph{minimal} if it is a shortest gallery with source $c_0'$ and target $c_{n+1}'$.
\end{definition}

We abbreviate ``$\gamma$ is a combinatorial gallery with source $x$ and target $y$'' by $\gamma:x\rightsquigarrow y$. In the following a \emph{gallery} is a combinatorial gallery in the sense of Definition~\ref{galleryDef}. 

We fix a type $t$ and denote by $\Ghat_t$ the set of all minimal galleries of type $t$ in $X$ with source $0\in A$. We denote by $\G_t$ the set of targets $\wgt(\gamma)$ of galleries $\gamma$  in $\Ghat_t$.
This leads to the following 
\begin{lemma}\label{preimage}
Let $A$ be an apartment with fixed special vertex $0$. Assume that $x$ is a special vertex in $A$ and $\gamma:0\rightsquigarrow x$ minimal of type $t$. Then
$$
r^{-1}(W.x) = \G_t.
$$
If $K\define\mathrm{Stab_{\mathrm{Aut(X)}}(0)}$ acts transitively on the set of all apartments containing $0$, then 
$$ r^{-1}(W.x) = K.x .$$
\end{lemma}

Note, for example, that $Stab_{\mathrm{Aut(X)}}(0)$ is transitive on the set of all apartments containing $0$ if $G$ has an affine $BN$-pair.

\begin{proof}
The retraction $r$ preserves adjacency relations and distance to $0$ in $X$. The first statement is an easy consequence of these properties. Secondly  transitivity of the $K$-action directly implies $\G_t=K.x$. 
\end{proof}

The following definition will help us describe the set $\rho(r^{-1}(W.x))$.
\begin{definition}\label{foldedGallery}
Fix an apartment $A$ in $X$ and let $H$ be a hyperplane, $d$ an alcove and $S$ a Weyl chamber in $A$. We say \emph{$H$ separates $d$ and $S$} if there exists a representative $S'$ of $\partial S$ in $A$ such that $S'$ and $d$ are contained in different half apartments determined by $H$.
\end{definition}

\begin{definition}
A gallery $\gamma= (c_0', c_0, c_1', c_1, c_2', \ldots , c_n', c_n, c_{n+1}' )$  is \emph{positively folded at $i$} if $c_i=c_{i-1}$ and 
the hyperplane $H\define \mathrm{span}({c_i'})$ separates $c_i$ and $\Cfm$.

A gallery $\gamma$ is \emph{positively folded} if it is positively folded at $i$ whenever $c_{i-1}=c_i$.
\end{definition}

The picture to have in mind is the following: Let $\gamma=(x, c_0,  c_1', c_1, \ldots , c_n, y )$ be a positively folded gallery of type $t$. Assume that there exists exactly one folding and let $i_0$ denote the associated index. Denote by $\widetilde{\gamma}=(x, d_0, d_1' , d_1,t \ldots , d_n, d_{n+1}' ) $ the minimal gallery of the same type $t$ such that $c_i=d_i$ and $c_i'=d_i'$ for all $i=0\ldots i_0-1$. Then $\gamma$ is obtained from $\widetilde{\gamma}$ by successively reflecting the remaining part $(d_{i_0}, d_{i_0+1}', \ldots , d_n, d_{n+1}')$ of $\widetilde{\gamma}$ at the hyperplane $H_{i_0}$.

\begin{lemma}\label{foldingLemma}
Given an element $\sigma$ of $\Ghat_t$. The image $\hat{\rho}(\sigma)$ is a positively folded gallery of the same type as $\sigma$.
\end{lemma}
\begin{proof}
Let $\gamma=(x, c_0, c_1', c_1, \ldots , c_n, y )$ be an image under $\hat{\rho}$ of an element $\sigma \in \Ghat_t$. Retractions preserve adjacency and hence the type of a gallery. Let $i$ be the smallest index such that $c_{i-1}=c_i$ and assume $\gamma$ is not positively folded at $i$. We may assume further that $\sigma$ has already been retracted onto $\gamma$ up to $i-1$, meaning 
$$
\sigma= (x, c_0 , c_1', c_1, \ldots , c_{i-1}, c_i', d_i, d_{i+1}', \ldots , d_n, y' ),
$$
where $d_j\notin A$ for all $j$. There is an apartment $A'$ containing $\Cfm, c_{i-1}$ and $d_i$. Denote by $H$ the hyperplane in $A'$ containing $c_i'$ and separating $d_i$ and $c_{i-1}$. One can find an alcove $b\in \Cfm\in A$ such that $\rho(c_i) = r_{A,b}(c_i)$. The retraction $r_{A,b}$ preserves distance to $b$ and maps $A'$ isometrically onto $A$. Hence the wall $H$ separates $c_i$ from $r_{b,A}(d_i)$ and also $r_{b,A}(d_i)$ from $\Cfm$. But this is not possible since we assumed that $r_{A,b}(d_i)=c_i=c_{i-1}$ and that $c_i$ is not separated from $\Cfm$ by $H$.
\end{proof}


\subsection*{Dual convexity}
Usually a subset $C$ of an apartment $A$ is defined to be $\sW$-convex if it is the intersection of finitely many half apartments. The \emph{$\sW$- convex hull} of a set $C$ is the intersection of all half apartments containing $C$. The following example illustrates why in our situation this is not the right notion of convexity. 

\begin{ex}
Let $A$ be a Coxeter complex of type $\widetilde{A}_2$. Let $\alpha_1$ and $\alpha_2$ be a basis of the underlying root system $R$. Let $x$ equal $3\alpha_1 + 3\alpha_2$. Observe that the $\sW$- convex hull of the Weyl group orbit $\sW.x$ contains $4\alpha_1+2\alpha_2$ which we denote by $y$. For special vertices $a,b$ let $\delta(a,b)$ be the length of the shortest gallery $\gamma$ such that $a$ is contained in the first and $b$ in the last chamber of $\gamma$. Then $\delta(0,x)=10$ and $\delta(0,y)=11$. 

For all elements $z$ in $\rho(r^{-1}(W.x))$ the distance $\delta(0,z)$ is less than or equal to $10$, since the two retractions are distance diminishing. Hence $y$ can not be contained in $\rho(r^{-1}(W.x))$.
\end{ex}

We therefore define a new notion of convexity using hyperplanes dual to the usual ones.

\begin{definition}
A \emph{dual hyperplane} is a set $H_{i,k}^*=\{x\in E \;\vert\; \langle \fw_{\alpha_i}, x^\vee \rangle=k\}$, where $k\in \Z$ and $i\in I\cup\{0\}$.
\end{definition}

Dual hyperplanes determine \emph{dual half apartments}. As $H_{i,k}$ is perpendicular to $\alpha$ the dual hyperplane $H_{i,k}^*$ is perpendicular to $\fw_i$. For any special vertex $x$ and root $\alpha\in R$ there exists a dual hyperplane perpendicular to $\alpha$ containing $x$. Notice that the positive cone is the intersection of all positve dual half-apartments $(H_{i,o}^*)^+$ where $\alpha$ is an element of $B$.

\begin{definition}
A \emph{convex} set is an intersection of finitely many dual half apartments in $A$, where the empty intersection is defined to be $A$. The \emph{convex hull} of a set $C$, denoted by $\dconv(C)$,  is the intersection of all dual half apartments containing $C$.
\end{definition}

Let $C$ be the orbit $W.x$ of a special vertex $x$ in $A$. One can prove that $\dconv(C)$ equals the metric convex hull of $W.x$ in $A$.

The following Lemma, which is a direct analogue of \cite[Lemma 3.3(2)]{Kostant}, plays a crucial role in the proof of our main result. Let $A^{type}(x)$ denote the set of all vertices in $\dconv(W.x)$ having the same type as $x$.

\begin{lemma}\label{conv^*}
Let $x$ be a special vertex in $A$ let $x^+$ denote the unique element of $W.x$ contained in $\Cf$. Then 
$$A^{type}(x)=\{y\in E : x^+-y^+\in ( \Cp\cap \mathcal{Q}) \}$$
where $\mathcal{Q}=\mathcal{Q}(R^\vee)$ is the co-weight-lattice of $R$. 
\end{lemma}
\begin{prf}
Assume without loss of generality that $x$ is contained in $\Cf$. Since $A^{type}(x)$ is $W$-invariant it suffices to prove that 
$$
A^{type}(x)\cap \Cf = \{y\in \Cf : x-y\in\Cp \text{ and } \type(x)=\type(y)\}.$$
Let $k_i=\langle \fw_i, x^\vee\rangle$, then $x=\sum k_i \fw_i$. Denote by $t$ the type of $x$. By the definition of dual hyperplanes,
\begin{equation}\label{num1}
A^{type}(x)\cap\;\Cf = \Cf \cap (\Cpm + x) = (\bigcap_{i=1}^n H_{i,0}^+) \cap (\bigcap_{i=1}^n (H^*_{i,k_i})^-).
\end{equation}

Let $y$ be an element of $A^{type}(x)\cap\;\Cf$. From (\ref{num1}) we deduce $\langle \fw_i, y^\vee\rangle \leq \langle \fw_i, x^\vee\rangle=k_i$ for all $i\in I.$ 
Consequently
\begin{equation*}
\langle \fw_i, (x-y)^\vee\rangle = \langle  \fw_i, x^\vee\rangle - \langle \fw_i, y^\vee\rangle \geq 0 
\Longleftrightarrow 
x-y \in \Cp.
\end{equation*}

Conversely assume that $y\in \Cf$ and $x-y\in \Cp$. Again one has
$$0\leq  \langle \fw_i, (x-y)^\vee \rangle
     = \langle \fw_i, x^\vee\rangle - \langle \fw_i, y^\vee\rangle
 $$
and hence $\langle \fw_i, x^\vee\rangle \geq \langle \fw_i, y^\vee\rangle$. Therefore $y\in A^{type}(x)$.
\end{prf}


\subsection*{A character formula}

The character formula stated below is used in the proof of our main result. For further references on representation theory see for example \cite{Bourbaki7-9} or \cite{Humphreys}.

Let $G$ be a connected semisimple complex algebraic group. Let $T$ be a maximal torus in $G$. Denote by $\mathfrak{X}$ the co-character group of $T$ and the subset of dominant co-characters in $\mathfrak{X}$ by $\mathfrak{X}_+$. Let $X$ be the affine Bruhat-Tits-building associated to $G(\C[\![t]\!])$ as defined in \cite{BruhatTits}. Apartments in $X$ can be identified with translates $gA$ of $A\define \mathfrak{X}\otimes_\Z\R$. Dominant co-weights $\lambda\in \mathfrak{X}_+$ correspond to special vertices in $\Cf$.
Denote by $V_\lambda$ the irreducible complex representation of highest weight $\lambda$ for $G$. Let $\chi_\lambda$ be the character of $V_\lambda$. 

In this setting Gaussent and Littelmann proved the following character formula. Note that our root system $R$ actually is the co-root system for the one considered in \cite{GaussentLittelmann}. Accordingly we switch co-characters to characters, etc.

\begin{thm}\cite{GaussentLittelmann}\label{characterFormula}
Fix $\lambda\in \mathfrak{X}_+$. Then
\begin{equation}
 \sum_{\gamma\in LS(t)} k^{\wgt(\gamma)} = \chi_\lambda = \sum_{\nu\in X} p_{\nu}k^\nu
\end{equation}
where $p_{\nu}$ is the dimension of the $v$-weightspace of $V_\lambda$ and $\wgt(\gamma)$ denotes the endpoint of $\gamma$.
\end{thm}

The set $LS(t)$ of $LS$-galleries, is a certain subset of the set of positively folded galleries of fixed type $t$ and source $0$ in $A$. The term $k^{wgt(\gamma)}$ in the formula is a 'group like element' of the group algebra associated with the targets $\wgt(\gamma)$ of $LS$-galleries. For all $\nu\in\Cf$ the coefficients $p_{\lambda}$ are strictly positive if and only if $\nu\leq \lambda$, which is equivalent to the fact that $\nu$ is contained in $A^{type}(x)\cap\Cf$. The fact that these group like elements are independent over the integers is crucial in the proof of Proposition~\ref{propA}. 

\begin{rem}
\cite{GaussentLittelmann} use a slightly more general notion of (combinatorial) galleries as given in Definition~\ref{galleryDef}. Nevertheless, by observations made in \cite{Schwer}, their result is true in our setting as well.
\end{rem}


\section{The main result}\label{convexityThm}

This section is devoted to the proof of our main result:

\begin{thm}\label{mainthm}
Let $X$ be a thick affine building with fixed apartment $A$. Let $x$ be a special vertex in $A$. Then
\begin{align*}
\rho(r^{-1}(\sW.x)) 	&\;=\dconv(W.x)\cap\{\text{ vertices } y\in A \;\text{of type } \mathrm{type}(x)\;\}\\
			&\;=\dconv(W.x)\cap (x+Q)\\
			&\;= A^\type(x).
\end{align*}
where $Q=Q(R^\vee)$ is the co-weight lattice of $R$.
\end{thm}

Extending the retractions $\rho$ and $r$ to minimal galleries in $X$ starting at $0$, the proof of Theorem~\ref{mainthm} can be reduced to the following two facts:

\begin{prop}\label{propA}
Let $A$ be an apartment of $X$ with fixed special vertex $0$ and fundamental Weyl chamber $\Cf$. Let $x$ be a special vertex in $A$ and let $x^+$ denote the unique element of $\sW.x$ contained in $\Cf$. Let $t$ be the type of a fixed minimal gallery $\gamma:0 \rightsquigarrow x^+$.
All vertices in $A^\type(x)=\{y\in \dconv(W.x)\cap (x+\mathcal{Q}(R^\vee)) \}$ are targets of positively folded galleries having type $t$. Conversely the targets of such galleries are contained in $A^\type(x)$. 
\end{prop}

\begin{prop}\label{propB}
Let $A$ be a fixed apartment of an affine building $X$. Fix an origin $0$ and fundamental Weyl chamber $\Cf$ in $A$. If $\gamma\subset A$ is a positively folded gallery with source $0$ of type $t$ then $\gamma$ has a pre-image under $\rho$ which is a minimal gallery starting at $0$. 
\end{prop}

With these propositions at hand, the proof of the main result reads as follows:

\begin{proofMainThm}
By Lemma~\ref{preimage}, the set $r^{-1}(W.x))$ equals $\G_t$. Therefore $\rho(r^{-1}(W.x)))=\rho(\G_t)$. Denote by $\hat{r}$ and $\hat{\rho}$ the extensions of $r$ and $\rho$ to minimal galleries starting at $0$. Lemma~\ref{foldingLemma} implies that $\hat{\rho}(\Ghat_t)$ is a set of positively folded galleries in $A$ and thus $\rho(\G_t)$ is the set of targets of the galleries contained in $\hat{\rho}(\Ghat_t)$. By Propositions~\ref{propA} and~\ref{propB} the assertion follows.
\end{proofMainThm}

\begin{rem}
Proposition~\ref{propA} is a purely combinatorial property of $LS$-galleries. During the preparation of the present paper Parkinson and Ram published their preprint \cite{ParkinsonRam}, where they give a proof the first assertion of Proposition~\ref{propA} using geometric and combinatorial properties of Coxeter complexes. The main idea of their proof carries over to thick non-discrete affine buildings. Details can be found in Section 8 of the authors thesis \cite{Diss}.
\end{rem}

It remains to prove Propositions~\ref{propA} and~\ref{propB}. 

\begin{proofPropB}
Let $\gamma=(0, c_0, d_1, c_1, d_2, \ldots ,d_n, c_n, y)$ be a positively folded gallery of type $t$ which is contained in $A$. We will inductively construct a pre-image of $\gamma$ under $\rho$ which is minimal in $X$.
Denote by $J\subset \{1,2,\ldots,n\}$ the set of indices $\{i_1, \ldots, i_k\}$ at which $\gamma$ is folded. Hence for each $i\in J$ we have that $c_{i_j}=c_{i_j -1}$. Assume $i_1 <i_2 <\ldots < i_k$. Let $H_{i_1}$ be the hyperplane in $A$ separating $c_{i_1}$ from $\Cfm$. We write $H_{i_1}^-$ (respectively $H_{i_1}^+$) for the half space of $A$ determined by $H_{i_1}$ containing a subsector of $\Cfm$ (respectively a subsector of $\Cf$). Since $X$ is a thick building there exists an apartment $A_1$ such that $A_1\cap A = H_{i_1}^-$. 

Define $A'\define H_{i_1}^+\cup(A_1\setminus A)$. 
By construction $\gamma$ is contained in $A'$. Let $r_{i_1}$ be the reflection associated to $H_{i_1}$ in $A'$ and define $\gamma^1$ as follows:
$$
\gamma^1=(0, c_0, d_1, \ldots  , d_{i_1}, c_{i_1}^1, d_{i_1+1}^1, \ldots , d_n^1, c_n^1, y^1)
$$
where $c_j^1 = r_{i_1}(c_j)$ and $d_j^1=r_{i_1}(d_j)$ for all $j\geq i_1$. Note that $d_{i_1}=d_{i_1}^1$ since $d_{i_1}$ is contained in $H_{i_1}$. Then $\gamma^1$ is a gallery of the same type as $\gamma$ which is not folded at $i_1$. Furthermore 
$\gamma^1$ is minimal up to the next folding index, i.e. the shortened gallery 
$(0, c_0, d_1, \ldots  , d_{i_1}, c_{i_1}^1, d_{i_1+1}^1, \ldots , d_{i_2 -1}^1, c_{i_2 -1}^1)$ 
is minimal. This step reduced the number of folding indices by one.
The construction also directly implies that $\rho(\gamma^1)=\gamma$.

We shorten the gallery $\gamma$ to 
$$
\gamma_1 = (x^1, c_{i_1}^1, d_{i_1+1}^1, \ldots , d_n^1, c_n^1, y^1)
$$
where $x^1$ is a vertex of $c_{i_1}$, which is chosen such that 
$(x^1, c_{i_1}^1, d_{i_1+1}^1, \ldots , d_{i_2}^1, c_{i_2-1}^1)$ 
is a minimal gallery from $x^1$ to $c_{i_2-1}^1$. By definition $\gamma^1$ is positively folded of suitably shortened type and contained in $A_1$.
As in the first step we can construct a gallery 
$$
\widetilde{\gamma}^2=(x_1, c_{i_1}^1, d_{i_1+1}^1, \ldots  , d_{i_2}^1, c_{i_2}^2, d_{i_2+1}^2, \ldots , d_n^2, c_n^2, y^2)
$$
 having analogous properties. Denote by $\rho_1$ the retraction $\rho_{A_1, \Cfm}$ onto $A_1$. The construction above implies that
$\rho_1(\widetilde{\gamma}^2)=\gamma^1$. 

Iterating the procedure $k$ times, we get a gallery
\begin{align*}
\gamma^k = 	&(0 , c_0 , d_1 
		 \ldots d_{i_1}, c_{i_1}^1, d_{i_1+1}^1 \ldots, 
		 \ldots d_{i_k}^{k-1}, c_{i_k}^k, d_{i_k+1}^k ,
		\ldots d_n^k, c_n^k, y^k)
\end{align*}
and a sequence of apartments $A_1, \ldots A_k$ such that $c_*^j, d_*^j\subset A_j$ for all $j=1,\ldots, k$. By construction $\gamma^k: 0\rightsquigarrow y^k$ is minimal in $X$ and of the same type as $\gamma$. We refer to the retraction $\rho_{A_{j}, \Cfm}$ onto the apartment $A_{j}$ as $\rho_j$. Lemma~\ref{compatibleRetractions} then implies that $\rho(\gamma^k)=(\rho_k\circ\rho_{k-1}\circ\ldots\circ\rho_1)(\gamma)$. Hence $\gamma^k$ is the desired preimage of $\gamma$.
\end{proofPropB}

\begin{lemma}\label{positiveFold}
Given an apartment $A$ with fixed origin an dfundamental Weyl chamber $\Cf$. Let $x$ be an element of $A$, let $H$ be a hyperplane separating $0$ and $\Cfm$ and denote by $H^+$ the half-apartment determined by $H$ which contains  a sub-Weyl chamber of $\Cfm$. Is $y$ a vertex in $\dconv(\sW.x)\cap H^+$ then the refected image of $y$ at $H$ is contained in $\dconv(\sW.x)$. 
\end{lemma}
\begin{proof}
Let $H'$ be the hyperplane parallel to $H$ containing $0$ and let $\alpha\in R$ be the root perpendicular to $H'$. The convex set $\dconv(\sW.x)$ is by definition $\sW$-invariant. Therefore the reflected image of $y$ at $H'$ is contained in $\dconv(sW.x)$ and equals $y-\langle y,\alpha^\vee\rangle \alpha$. Since $H$ separates $0$ and $y$ the reflected image of $y$ at $H$ equals $y-\lambda \alpha$ for some $\lambda < \langle y,\alpha^\vee\rangle$
and is therefore containe in $\dconv(\sW.x)$ as well.
\end{proof}

\begin{proofPropA}
Let $\Gamma_t^+$ be the set of all positively folded galleries of type $t$ in $A$ starting in $0$. This set contains $LS(t)$, the so called $LS$-galleries (see \cite{GaussentLittelmann} for the definition). Let $G$ be a connected semisimple complex algebraic group $G$, denote by $\mathfrak{X}$ the character group of a chosen maximal torus $T$. Let $G$ be chosen such that $A$ can be identified with $\mathfrak{X}\otimes_\Z\R$. 
By Theorem~\ref{characterFormula} we have the following formula for $\chi_{x^+}$, the character of the (irreducible) highest weight representation $V_{x^+}$ with highest weight $x^+$:
\begin{equation*}
\sum_{\sigma\in LS(t)} k^{\wgt(\sigma)} = \chi_{x^+} = \sum_{\nu\in X} p_{\nu}k^\nu
\end{equation*}
where $\wgt(\sigma)$ is the target of $\sigma$. The coefficients $p_{\nu}$ in the right part are the dimensions of the $\nu$-weight spaces of $V_{x^+}$, which are positive if and only if $\nu\leq x^+$. Simply comparing coefficients of $k^\nu$ on the left and right hand side of the equation, there has to exist an $LS$-gallery with target $\nu$ for all $\nu\leq x^+$. By Lemma~\ref{conv^*} the proposition follows.

First observe that Theorem 3 and Lemma 11 of \cite{GaussentLittelmann} imply that $\rho$ maps a minimal gallery $\gamma$ of fixed type with source $0$ onto a positively folded gallery of the same type. Recall that $\gamma$ is the pre-image under $r$ of a minimal gallery $\sigma$ (of the same type) from $0$ to some point in $\sW.x$. Choosing a preimage $\gamma$ of $\sigma$ and applying $\rho$ is the same as folding $\sigma$ at certain panels contained in hypeplanes separating $0$ and $\Cfm$, i.e. folding once at panel $d$ corresponds to a reflection along a hyperplane $H$ containing $d$. The endpoint $y$ of the unfolded gallery is mapped onto $r(y)$, where $r$ denotes the reflection at $H$. Therefore the endpoint of a folded gallery is the reflected image of an orbit point $y^0\in \sW.x$ under a sequence of reflections $r^i$ along hyperplanes $H^i$ all separating $0$ and the current point $y^{i-1}$. Unfolding the gallery far enough, similar to the argument in the proof of \ref{propB}, we may assume that there exists one folding index only. Then Lemma \ref{compatibleRetractions} and  Lemma \ref{positiveFold} imply the assertion.
\end{proofPropA}


\section{An application}

Let $G$ be a group with affine $BN$-pair. There exists then an affine building $(X,\App)$ associated to $G$ on which the group acts. Let $B$ and $N$ denote the groups of the sphercial $BN$-pair of $G$ which we assume to be \emph{split}, i.e. $B=UT$. Denote by $\Delta$ the spherical building associated to $G$ by means of the spherical $BN$-pair. Notice that $\Delta$ is the building at infinity of $X$. The group $N$ stabilizes an apartment $a$ of $\Delta$ and $B$ stabilizes a chamber $c$ of $a$.  The groups $U$ and $T$ can be interpreted as follows: $T$ is a group of translations in the affine apartment $A$ with boundary $\partial A = a$ and $U$ acts simply transitive on the set of all affine apartments containing $c$ at infinity. Fix a chart of $A$ such that $c=\partial(\Cfm)$.

Let $K$ be the stabilizer of $0\in A$ in $G$. Then $G$ has an \emph{Iwasawa decomposition}\index{Iwasawa decomposition}
$$ G=BK=UTK .$$ 
Notice that special vertices of the same type as $0$ in $X$ are in one-to-one correspondence with cosets of $K$ in $G$. Let $x$ be a special vertex in $X$ having the same type as $0$.
Then there exists an element $u\in U$ such that $(u^{-1}).x\in A$. Let $t\in T$ be the translation mapping $0$ to $(u^{-1}).x$. If the origin $0$ is identified with $K$ the vertex $x$ corresponds to the coset $utK$. Denote by $t_0$ the type of the origin. One has :
\begin{align*}
\{ \text{special vertices of type } t_0 \text{ in } X\} &\stackrel{1:1}{\longleftrightarrow} \{ \text{ cosets of } K \text{ in } G \}\\
 0 	& \longmapsto  K \\
X\ni x 	& \longmapsto  utK  \text{ with } u,t \text{ chosen as above.}
\end{align*}
Any special vertex of type $t_0$ of $X$ can hence be identified with a coset $utK$ with suitably chosen $u\in U$ and $t\in T$. Note that vertices of type $t_0$ contained in $A$ correspond precisely to cosets of the form $tK$ with $t\in T$. 

Furthermore it is easy to see that $\rho$ is exactly the projection that maps $utK$ to $tK$, and hence to $A$. 

For all special vertices $tK$ contained in $A$ the set $r^{-1}(\sW.tK)$ is the same as the $K$-orbit of $tK$ which is 
$$r^{-1}(\sW.tK)= KtK.$$ 

In the situation as described above the following theorem is a direct reformulation of Theorem~\ref{mainthm}.

\begin{thm}\label{Thm_BNpair}
For all $tK\in A$ we have 
$$\rho(KtK)=\dconv(\sW.tK)$$
or, since $\rho^{-1}(t'K)=Ut'K$, equivalently
$$
\emptyset \neq Ut'K \cap KtK \:\Longleftrightarrow\: t'K \in \dconv(\sW.tK).
$$
\end{thm}

The proof of ``$\Rightarrow$'' in the second statement is well known and can be found in \cite{BruhatTits}. As already mentioned in the introduction partial results on ``$\Leftarrow$'' and related questions can, for example, be found in \cite{GaussentLittelmann, MirkovicVilonen, Rapoport, Schwer} or \cite{Silberger}.

\address{Petra Hitzelberger, Mathematisches Institut der Westf\"alischen Wilhelms-Universit\"at M\"unster, Einsteinstrasse 62, 48149 M\"unster.}
\email{hitzelberger@uni-muenster.de}
~\clearpage

\end{document}